\documentclass{article}

\usepackage{amssymb}
\usepackage{amsthm}

\newtheorem{thm}{Theorem}[section]
\newtheorem{prop}[thm]{Proposition}
\newtheorem{lem}[thm]{Lemma}
\newtheorem{cor}[thm]{Corollary}

\newtheorem{defn}{Definition}

\newcommand{\sA}{\mathcal{A}}
\newcommand{\rat}{\mathbb{Q}}

\begin{document}

\title{Classification From a Computable Viewpoint}

\author{W.\ Calvert and J.\ F.\ Knight}

\maketitle

\section{Introduction}

Classification is an important goal in many branches of mathematics.
The idea is to describe the members of some class of mathematical
objects, up to isomorphism or other important equivalence, in terms of
relatively simple invariants.  Where this is impossible, it is useful
to have concrete results saying so.  In model theory and descriptive
set theory, there is a large body of work, showing that certain
classes of mathematical structures admit classification, while others
do not.  In the present paper, we describe some recent work on
classification in computable structure theory.

Section 1 gives some background from model theory and descriptive set
theory.  From model theory, we give sample structure and non-structure
theorems for classes that include structures of arbitrary cardinality.
We also describe the notion of Scott rank, which is useful in the more
restricted setting of countable structures.  From descriptive set
theory, we describe the basic Polish space of structures for a fixed
countable language with fixed countable universe.  We give sample
structure and non-structure theorems based on the complexity of the
isomorphism relation, and on Borel embeddings.

Section 2 gives some  background on 
computable structures.  We
describe three approaches to classification for these structures.  The
approaches are all equivalent.  However, one approach, which involves
calculating the complexity of the isomorphism relation, has turned out
to be more productive than the others.  Section 3 describes results on
the isomorphism relation for a number of mathematically interesting
classes---various kinds of groups and fields.  In Section 4, we
consider a setting similar to that in descriptive set theory.  We
describe an effective analogue of Borel embedding which allows us to
make distinctions even among classes of finite structures.  Section 5
gives results on computable structures of high Scott rank.  Some
of these results make use of computable embeddings.

We shall assume some background in computability---$m$-reducibility, 
arithmetical and hyperarithmetical hierarchies, $\Sigma^1_1$ and 
$\Pi^1_1$ sets and relations.  This material may be found in
\cite{R}, or in \cite{A-K}.  

\subsection{Background from model theory}             

In classical model theory, the basic setting is the class $K$ of models for a
countable complete elementary first order theory.  If the theory has
an infinite model, then there are models of all infinite
cardinalities.  For the theory of vector spaces over the rationals,
each model is determined up to isomorphism by its dimension.  More
generally, for an $\aleph_1$-categorical theory $T$, each model is
determined, up to isomorphism, by the dimension of a ``strongly
minimal'' formula, with parameters satisfying a principal type
\cite{Mo1}, \cite{B-L}.  This is the prototypical structure theorem.
For the theory of dense linear orderings without endpoints, or other
unstable theory, there are $2^\kappa$ non-isomorphic models of
cardinality $\kappa$, for all $\kappa\geq \aleph_1$, too many to allow
nice invariants \cite{Sh1}.  (See \cite{Ho} for a discussion of
further structure and non-structure theorems in the setting of models
having arbitrary cardinality.)

Vaught \cite{V} focused attention on countable models.  Let $T$ be a countable complete elementary first order theory, and consider the number of countable models, up to isomorphism.  There are familiar examples illustrating some possible numbers.  For the theory of dense linear orderings without endpoints, the number is $1$.  For the theory of vector spaces over the rationals, the number is $\aleph_0$.  For true arithmetic, the number is $2^{\aleph_0}$.  Ehrenfeucht found an example for which the number is $3$:  the theory of dense linear orderings without endpoints, with a strictly increasing sequence of constants. The example is easily modified to give examples for which the number is $n$, for all $n\geq 3$.  The Ehrenfeucht examples are described in \cite{V}.  Vaught showed that for a countable complete theory, the number of isomorphism types of countable models cannot be $2$.  

Vaught conjectured that for a countable complete theory, the number of countable models, up to isomorphism is either $\leq\aleph_0$ or $2^{\aleph_0}$.  Vaught's Conjecture has been proved for many special kinds of theories, including theories of linear orderings \cite{Ru} and trees \cite{St}, $\omega$-stable theories \cite{S-H-M}, and superstable theories of finite rank \cite{Bu}.   
In the setting of Vaught's Conjecture, the statement that $T$ has $2^{\aleph_0}$ isomorphism types of countable models is a non-structure theorem.  The statement that $T$ has only countably many isomorphism types of countable models is a kind of structure theorem.  While it does not provide descriptions for the models, it holds out the possibility that there may be nice descriptions.

\subsection{Infinitary formulas}       

In describing countable structures, certain infinitary (but still first order) sentences are useful.  For a language $L$, the \emph{$L_{\omega_1\omega}$} formulas are infinitary formulas in which the disjunctions and conjunctions over countable sets.  For a thorough discussion of $L_{\omega_1\omega}$, see Keisler's beautiful book \cite{Ke}.  Scott \cite{Sc} proved the following.  

\begin{thm} [Scott Isomorphism Theorem]
\label{thm1.1}

Let $\mathcal{A}$ be a countable structure (for a countable language $L$).  Then there is an $L_{\omega_1\omega}$ sentence $\sigma$ whose countable models are just the isomorphic copies of $\mathcal{A}$.  

\end{thm}

A sentence $\sigma$ with the property of Theorem \ref{thm1.1} is called a \emph{Scott sentence} for $\mathcal{A}$.  In proving Theorem \ref{thm1.1}, Scott assigned countable ordinals to tuples in $\mathcal{A}$, and to $\mathcal{A}$ itself.  There are several different definitions of \emph{Scott rank} in use.  We begin with a family of equivalence relations.  

\begin{defn}  

Let $\overline{a}$, $\overline{b}$ be tuples in $\mathcal{A}$.

\begin{enumerate}

\item  $\overline{a}\equiv^0\overline{b}$ if $\overline{a}$ and $\overline{b}$ satisfy 
the same quantifier-free formulas,

\item  for $\alpha > 0$, $\overline{a}\equiv^\alpha\overline{b}$ if
  for all $\beta < \alpha$, for each $\overline{c}$, there exists
  $\overline{d}$, and for each $\overline{d}$, there exists
  $\overline{c}$, such that
  $\overline{a},\overline{c}\equiv^\beta\overline{b},\overline{d}$.

\end{enumerate}

\end{defn}

Our relations $\equiv^\alpha$ differ from Scott's in that we consider
tuples $\overline{c}$ and $\overline{d}$, where Scott considered
only single elements $c$ and $d$.

\begin{defn}\

\begin{enumerate}

\item  The \emph{Scott rank} of a tuple $\overline{a}$ in $\mathcal{A}$ is the least $\beta$ such that for all $\overline{b}$, $\overline{a}\equiv^\beta\overline{b}$ implies $(\mathcal{A},\overline{a})\cong (\mathcal{A},\overline{b})$. 

\item  The \emph{Scott rank} of the structure $\mathcal{A}$, denoted by $SR(\mathcal{A})$, is the least ordinal $\alpha$ greater than the ranks of all tuples in $\mathcal{A}$.

\end{enumerate}

\end{defn}

\noindent
\textbf{Example}.  If $\mathcal{A}$ is an ordering of type $\omega$, then $SR(\mathcal{A}) = 2$.

\bigskip

Morley \cite{Mo2} showed that if $T$ is a counter-example to Vaught's Conjecture, then the number of isomorphism types of countable models must be $\aleph_1$.  The idea is behind the proof is that each model must have countable Scott rank, and for each countable ordinal $\alpha$, the number of models of rank $\alpha$ must be countable or $2^{\aleph_0}$.  Vaught's Conjecture holds for $T$ iff for some countable ordinal $\alpha$ either there are $2^{\aleph_0}$ non-isomorphic models of rank $\alpha$, or else all countable models have Scott rank at most $\alpha$.  Shelah's proof of Vaught's Conjecture for $\omega$-stable theories \cite{S-H-M} involves putting a fixed countable ordinal bound on the Scott ranks of all models.  For more information related to Scott rank, see \cite{Sa}.

\subsection{Background from descriptive set theory}

Next, we turn to descriptive set theory.  We mention only a few notions and results.  For more information, see \cite{Hj-K} or \cite{B-K}.  We describe a version of the basic setting.    
The structures are countable, with fixed universe $\omega$.  The languages are also countable.  
For a given language $L$, we have a list $(\varphi_n)_{n\in\omega}$ of all atomic sentences involving symbols from $L$ plus constants from $\omega$.  We may identify a
structure $\mathcal{A}$ with the function $f\in 2^\omega$ such that 
\[f(n) = \left\{\begin{array}{ll}
1 & \mbox{ if $\mathcal{A}\models \varphi_n$}\\
0 & \mbox{ if $\mathcal{A}\models\neg{\varphi_n}$}
\end{array}
\right.\]  

The class of all structures for the language, or functions in
$2^\omega$, has a natural metric topology, where the distance between
$f,g\in 2^\omega$ is $1/2^n$ if $n$ is least such that $f(n)\not=
g(n)$.   The class of structures for the language forms a \emph{Polish
space}; i.e., it is separable and complete.  The  basic open
neighborhoods have the form 
$N_\sigma = \{f\in 2^\omega:f\supseteq\sigma\}$, where 
$\sigma\in 2^{<\omega}$.  The
\emph{Borel} subsets are generated from the basic open neighborhoods,
by taking countable unions and complements.  We may also consider Borel
subsets of a product of classes.  We have Borel relations on a given class, and Borel functions from one class to another.

We consider classes $K$, contained in the Polish space described above, such that $K$ is closed under isomorphism.  By results of Scott \cite{Sc}, Vaught \cite{V2}, and D.\ Miller \cite{DM}, $K$ is Borel iff it is the class of models of some $L_{\omega_1\omega}$ sentence.  We add the assumption that $K$ is Borel.  One approach to classification in this setting is to consider the isomorphism relation on $K$.  A result saying that the isomorphism relation on a class $K$ is
Borel is a weak kind of structure theorem, and a result saying that the isomorphism
relation is not Borel is a non-structure theorem.    

For the class of vector spaces over $\mathbb{Q}$, or the class of algebraically closed fields of a given characteristic, the isomorphism relation is Borel.  More generally, for any class $K$ with only countably many isomorphism types, the isomorphism relation is Borel.  The converse is not true.  For example, if $K$ is the class of Archimedean ordered fields, then there are $2^{\aleph_0}$ different isomorphism types, but the isomorphism relation is still Borel.  For the class of linear orderings, the isomorphism relation is not Borel.  Perhaps surprisingly, Hjorth \cite{Hj} showed that for the class of torsion-free Abelian groups, the isomorphism relation is not Borel.

There is a second approach to classification, which involves comparing classes, and provides a great deal of information.  This is the approach of Friedman and Stanley \cite{F-S}.

\begin{defn}

Let $K,K'$ be classes of structures.  We say that $\Phi$ is a \emph{Borel embedding} of $K$ in $K'$ if $\Phi$ is a Borel function from $K$ to $K'$, and for $\mathcal{A},\mathcal{A}'\in K$, $\mathcal{A}\cong\mathcal{A}'$ iff $\Phi(\mathcal{A})\cong\Phi(\mathcal{A}')$.  (More properly, we might say that $\Phi$ is an embedding of $K$ in $K'$ \emph{up to isomorphism}, or an embedding of $K/_{\cong}$ in $K'/_{\cong}$.)  

\end{defn}

\noindent
\textbf{Notation}:  We write $K\leq_B K'$ if there is a Borel embedding of $K$ in $K'$.  

\begin{prop}
\label{prop1.2}

The relation $\leq_B$ is a partial order on the set of all classes. 

\end{prop} 
    
A Borel embedding gives us a way of transferring descriptions from one class to another.  If $K\leq_B K'$, under the Borel embedding $\Phi$, and $K'$ has simple invariants, then we may describe $\mathcal{A}\in K$ by giving the invariants for $\Phi(\mathcal{A})$.  
We say that $K$ is \emph{Borel complete} if for all $K'$, $K'\leq_B K$; equivalently, $K'\leq_B K$, where $K'$ is the class of all structures for a language with at least one at least binary relation symbol.  Our intuition says that $K'$ cannot be classified, so the statement that $K$ is Borel complete is a convincing non-structure theorem.  Friedman and Stanley \cite{F-S} produced Borel embeddings of undirected graphs in fields of arbitrary characteristic, trees, and linear orderings, showing that all of these classes are ``Borel complete''.  We shall say more later about some of these embeddings.              

\section{Computable structures}  

We now turn to computable structures.  The languages that we consider are computable; i.e., the set of non-logical symbols is computable, and we can effectively determine the kind (relation or function) and the arity.  The structures have universe a subset of $\omega$, and we identify a structure $\mathcal{A}$ with its atomic diagram $D(\mathcal{A})$.  A structure $\mathcal{A}$ is \emph{computable} if $D(\mathcal{A})$ is computable. 
Then a \emph{computable index} for $\mathcal{A}$ is a number $e$ such that
$\varphi_{e}$ is the characteristic function of $D(\mathcal{A})$.  We write $\mathcal{A}_{e}$ for the structure with index $e$ (if there is such a structure).                      

\subsection{Computable infinitary formulas}

In describing computable structures, it is helpful to use \emph{computable infinitary formulas}.  Roughly speaking, these are infinitary formulas in which the disjunctions and conjunctions are over computably enumerable sets.  For a more thorough discussion, see~\cite{A-K}. 

\bigskip
\noindent
\textbf{Example}: There is a natural computable infinitary sentence whose
models are just the Abelian $p$-groups---the conjunction of the
usual axioms for Abelian groups, plus
$$\forall x\,\bigvee_{n}\hspace{-.22in}\bigvee\ 
\underbrace{x+\ldots+x}_{p^{n}}\  =\  0\ .$$   

\bigskip

All together, computable infinitary formulas have the same
expressive power as those in the least admissible fragment of
$L_{\omega_{1},\omega}$.  The usefulness comes from the following classification.

\bigskip
\noindent
\textbf{Complexity of formulas}

\begin{enumerate}

\item A computable $\Sigma_{0}$ and $\Pi_{0}$ formula is a finitary open formula.

\item  For $\alpha > 0$, a computable $\Sigma_{\alpha}$ formula $\varphi(\overline{x})$ is a c.e.\ disjunction of formulas of the form $\exists\overline{u}\,\psi(\overline{u},\overline{x})$, where
$\psi$ is computable $\Pi_{\beta}$ for
some $\beta<\alpha$.

\item  For $\alpha > 0$, a computable $\Pi_{\alpha}$ formula is a c.e.\ conjunction of formulas of the form $\forall\overline{u}\,\psi(\overline{u},\overline{x})$, where 
$\psi$ is computable $\Sigma_\beta$ for
some $\beta<\alpha$.  
 
\end{enumerate}

\begin{prop}
\label{prop2.1}

Satisfaction of computable $\Sigma_{\alpha}$
(or $\Pi_{\alpha}$) formulas in computable structures is $\Sigma^{0}_{\alpha}$ (or $\Pi^{0}_{\alpha}$), with all possible uniformity. 

\end{prop}

There is a version of Compactness for computable infinitary formulas which, unlike the usual Compactness, can be used to produce computable structures.
The original result of this kind, for weak second order logic, is due to Kreisel (in a footnote in \cite{Kr}).  Barwise \cite{B} gave a result for countable admissible fragments of $L_{\omega_1\omega}$.  The special case below, with its corollaries, may be found in \cite{A-K}.   

\begin{thm} [Barwise-Kreisel Compactness] 
\label{thm2.2}

If $\Gamma$ is a $\Pi^{1}_{1}$
set of computable infinitary sentences and every $\Delta^{1}_{1}$ subset
has a model, then $\Gamma$ has a model.

\end{thm}

\begin{cor}
\label{cor2.3}

If $\Gamma$ is a $\Pi^{1}_{1}$
set of computable infinitary sentences and every $\Delta^{1}_{1}$ subset
has a computable model, then $\Gamma$ has a computable model.

\end{cor}

The next two corollaries illustrate the power of computable infinitary formulas to describe computable structures.

\begin{cor}
\label{cor2.4}

If $\mathcal{A}$ and $\mathcal{B}$ are computable (or hyperarithmetical) structures satisfying the same
computable infinitary sentences, then $\mathcal{A}\cong\mathcal{B}$. 

\end{cor}

\begin{cor}[Nadel \cite{N1}]
\label{cor2.5}

Let $\mathcal{A}$ be a computable (or hyperarithmetical) structure.  If $\overline{a}$ and $\overline{b}$ are tuples in $\mathcal{A}$ satisfying the same computable infinitary formulas, then there is an automorphism of $\mathcal{A}$ taking $\overline{a}$ to $\overline{b}$. 

\end{cor}

The final corollary is a special case of a stronger result, due to Ressayre \cite{Re}.

\begin{cor} [Ressayre]
\label{cor2.6}

Let $\mathcal{A}$ be a computable (or hyperarithmetical) structure,
and let $\Gamma$ be a $\Pi^1_1$ set of computable infinitary sentences
involving finitely many new symbols, in addition to those in the
language of $\mathcal{A}$.  If every $\Delta^1_1$ subset of $\Gamma$
is satisfied in an expansion of $\mathcal{A}$, then $\mathcal{A}$ has
an expansion satisfying $\Gamma$.

\end{cor}

\subsection{Structure and non-structure theorems}

Certain classes of structures arising in computable structure theory
can be classified simply.  Goncharov \cite{G1}, \cite{G2} showed that
the Boolean algebras with a property of ``computable categoricity''
are the ones with only finitely many atoms, and the linear orderings
with this property are the ones with finitely many successor pairs.
In a talk in Kazan in 1997, Goncharov stated a large number of further
problems, calling for classification of further classes of computable
structures.  Some of these problems seemed likely to have nice
answers, while others did not.  At the end of the talk, Shore asked
Goncharov what would make him give up.  Shore's question was
considered in \cite{G-K}.  The aim was to find convincing statements
for structure and non-structure theorems, like those in model theory
and descriptive set theory.

Let $K$ be a class of structures, closed under isomorphism.  Let
$K^{c}$ be the set of computable elements of $K$.  Let $I(K)$ be the
set of computable indices for elements of $K^{c}$.  We suppose that
there is a computable infinitary sentence for which $K^c$ is the class
of computable models.  Equivalently, $I(K)$ is hyperarithmetical.  In
\cite{G-K}, three different approaches are discussed.  All give
intuitively correct answers for some familiar classes (e.g.,
$\mathbb{Q}$-vector spaces can be classified, linear orderings
cannot).

The first approach involves ``Friedberg enumerations''.  

\begin{defn}

An \emph{enumeration} of $K^c/_{\cong}$ is a sequence $(\mathcal{C}_n)_{n\in\omega}$ of elements of $K^c$ such that for all 
$\mathcal{A}\in K^c$, there exists $n$ such that $\mathcal{C}_n\cong\mathcal{A}$; i.e., the sequence represents all isomorphism types of computable members of $K$.  The enumeration is \emph{Friedberg} if for each $\mathcal{A}\in K^c$, there is a unique $n$ such that $\mathcal{C}_n\cong\mathcal{A}$; i.e., each isomorphism type appears just once.  

\end{defn}

The complexity of the enumeration is that of the sequence of
computable indices for elements of $K$.  Saying that $K^c/_{\cong}$
has a hyperarithmetical Friedberg enumeration is an assertion of
classifiability.  Lack of repetition in the list reflects
understanding of the different isomorphism types.  This approach has
considerable appeal.  For one thing, in mathematical practice, a
classification is often given by a list---think of the classification
of finite simple groups.  Another attractive feature of this approach
is that it fits with an established body of work on enumerations, in
computability and computable structure theory \cite{E}.

The second approach to classification considered in \cite{G-K} is
related to Scott ranks and Scott sentences.  The result below, a
consequence of Corollary \ref{cor2.5}, says which Scott ranks are
possible for computable structures.  We have already mentioned
computable ordinals.  These are the order types of computable well
orderings of $\omega$, and they form a countable initial segment of
the ordinals.  The first non-computable ordinal is denoted by
$\omega_1^{CK}$.  Nadel used a different definition of Scott rank, but
his work in \cite{N1} shows the following.

\begin{thm}
\label{thm2.6}

For a computable structure $\mathcal{A}$,
$SR(\mathcal{A})\leq\omega_1^{CK}+1$. 

\end{thm}

There are computable structures taking the maximum possible rank.
Harrison \cite{Ha} proved the following.

\begin{prop} [Harrison]
 \label{harrisonorder}

There is a computable ordering of type\linebreak 
$\omega_1^{CK}(1+\eta)$ (with a first interval of type $\omega_1^{CK}$, followed by densely many more).  

\end{prop}

\begin{proof}

Kleene showed that there is a computable tree $T\subseteq\omega^{<\omega}$ with a path but no hyperarithmetical path.  Let $<$ be the Kleene-Brouwer ordering on $T$.  Under this ordering, $\sigma <\tau$ iff either $\sigma$ properly extends $\tau$ or else there is some first $n$ on which $\sigma$ and $\tau$ differ, and for this $n$,
$\sigma(n) < \tau(n)$.  Harrison showed that $(T,<)$ has order type
$\omega_1^{CK}(1+\eta) + \alpha$, where $\alpha$ is a computable
ordinal.  There is an initial segment of type
$\omega_1^{CK}(1+\eta)$.

\end{proof} 

The Harrison ordering has Scott rank $\omega_1^{CK}+1$.    
If $a$ is an element outside the initial interval of type $\omega_1^{CK}$, then $a$ has Scott rank $\omega_1^{CK}$. 

\begin{thm}

For a computable structure $\mathcal{A}$,

\begin{enumerate}

\item  $SR(\mathcal{A}) < \omega_1^{CK}$ if there is some computable ordinal $\beta$ such that the orbits of all tuples are defined by computable $\Pi_\beta$ formulas.

\item  $SR(\mathcal{A}) = \omega_1^{CK}$ if the orbits of all tuples are defined by computable infinitary formulas, but there is no bound on the complexity of these formulas.

\item  $SR(\mathcal{A}) = \omega_1^{CK}+1$ if there is some tuple whose orbit is not defined by any computable infinitary formula.

\end{enumerate}

\end{thm} 

If a computable structure has low Scott rank, then we expect it to have a simple Scott sentence.  Nadel \cite{N1}, \cite{N2} showed the following.

\begin{thm} [Nadel]
\label{thm2.7}

A computable (or hyperarithmetical) structure has computable Scott rank iff it has a computable infinitary Scott sentence.

\end{thm}

\begin{prop}
\label{prop2.8}

Suppose $K$ is a class of structures closed under isomorphism.
If there is a computable bound on the Scott ranks of elements of $K^c$, then there is a computable ordinal $\alpha$ such that the elements of $K^c$ all have computable $\Pi_\alpha$ Scott sentences.     

\end{prop}

If we understand that goal of computable classification is to describe each computable member of $K$, so as to distinguish its computable copies from other \emph{computable} structures, then something weaker than a Scott sentence may suffice.  A \emph{pseudo-Scott sentence for $\mathcal{A}$} is a sentence whose computable models are just the computable copies of $\mathcal{A}$.  We do not know whether there is a computable structure with a computable infinitary pseudo-Scott sentence but no computable infinitary Scott sentence.  What we can prove is the following.   

\begin{prop}
\label{prop2.9}

Let $K$ be a class of structures closed under isomorphism.  Then the following are equivalent.

\begin{enumerate}

\item  There is a computable ordinal $\alpha$ such that any elements of $K^c$ satisfying the same computable $\Pi_\alpha$ sentences are isomorphic.  

\item  Each $\mathcal{A}\in K^c$ has a computable infinitary pseudo-Scott sentence.      

\end{enumerate}  

\end{prop}

The statements in Proposition \ref{prop2.9} are structure theorems.               

\bigskip

We turn to the third approach.      

\begin{defn}

The \emph{isomorphism problem} for $K$ is 
$$E(K) = \{(a,b):a,b\in I(K)\ \&\ \mathcal{A}_{a}\cong\mathcal{A}_{b}\}\ .$$

\end{defn}

The statement that $E(K)$ is hyperarithmetical is one more structure theorem.  For the class of computable models of a computable infinitary sentence, the 
three approaches are equivalent.

\begin{thm} 
\label{thm2.9} 

If $I(K)$ is hyperarithmetical, then the following are equivalent:

\begin{enumerate}

\item  there is a computable ordinal $\alpha$ such that any two members of $K^c$ satisfying the same computable $\Pi_\alpha$ sentences are isomorphic,

\item  $E(K)$ is hyperarithmetical,

\item  $K^c/_{\cong}$ has a hyperarithmetical Friedberg enumeration.

\end{enumerate} 

\end{thm}

\begin{proof}

It is not difficult to show that 1 implies 2 and 2 implies 3.  To see
that 3 implies 1, we form a hyperarithmetical structure
$\mathcal{A}^*$ including all structures in $K^c$, made disjoint.  We
also include a set corresponding to $I(K)$ and a relation $\mathbb{Q}$
associating each index $a\in I(K)$ with the elements of the
corresponding structure $\mathcal{A}_a$.  Let $I$ be a
hyperarithmetical subset of $I(K)$, representing indices in the
Friedberg enumeration.  We add this to $\mathcal{A}^*$.  Let $a,b$ be 
new constants, and let
$\Gamma(a,b)$ be a $\Pi^1_1$ set of computable infinitary formulas
saying $a,b\in I$, $x\not= y$, and for all computable ordinals
$\alpha$, $\mathcal{A}_a$ and $\mathcal{A}_b$ satisfy the same
computable $\Pi_\alpha$ sentences.  If 1 fails, then every
$\Delta^1_1$ subset of $\Gamma(a,b)$ is satisfiable in
$\mathcal{A}^*$.  Then using Corollary \ref{cor2.6}, we can show
that the whole set is satisfiable.  This is a contradiction.

\end{proof}            
  
\section{Isomorphism problems}

Proposition \ref{prop2.9} says that the three approaches to
classification for computable structures are all equivalent.  The
first approach, involving Friedberg enumerations, is certainly
appealing.  However, the third approach, involving isomorphism
problems, has proved to be more productive.  When we attempt to 
determine, for various specific classes $K$, whether there is a 
hyperarithmetical Friedberg enumeration of $K^c/_{\cong}$, it seems that
we must first calculate the complexity of $E(K)$.  In this section, we give a 
number of these calculations, mainly from \cite{C1}, \cite{C2}.  

Of course, if $I(K)$ is hyperarithmetical, then $E(K)$ must be 
$\Sigma^1_1$, since the isomorphism relation may be defined by a 
statement that asserts the existence of a function that is total and bijective and
preserves the basic relations and operations.  For some classes of structures, it is
well-known that the isomorphism problem is $m$-complete $\Sigma^1_1$.    

\begin{thm} [Folklore] 
\label{thm3.1} 

For each of following classes $K$, $E(K)$
is $m$-complete $\Sigma^{1}_{1}$:

\begin{enumerate}

\item  linear orderings

\item undirected graphs

\item  Boolean algebras

\item  Abelian $p$-groups

\item  arbitrary structures for language with at least one at least
binary relation symbol.

\end{enumerate}

\end{thm}

The ideas go back to Kleene and the fact that the class of indices for computable trees without paths, and the class of indices for well orderings are $m$-complete $\Pi^1_1$ (see \cite{R}).  For 1, we show that for any $\Sigma^1_1$ set $S$, there is a uniformly computable sequence $(\mathcal{A}_n)_{n\in\omega}$ such that if $n\in S$, then $\mathcal{A}_n$ is a Harrison ordering, and if $n\notin S$, then $\mathcal{A}_n$ is a well ordering.  We have an $m$-reduction of $S$ to the set of indices for the Harrison ordering.  For parts 2-5, the proof is similar.  In each case, the class $K$ contains a structure $\mathcal{A}$ such that for any $\Sigma^1_1$ set $S$, there is a uniformly computable sequence $(\mathcal{A}_n)_{n\in\omega}$ in $K$ such that $\mathcal{A}_n\cong\mathcal{A}$ iff $n\in S$.  For complete proofs of all parts of the theorem, see \cite{G-K}.

\bigskip

Items 1, 2, and 5 in Theorem \ref{thm3.1} match our intuition.  Whatever classification means, it should be impossible for these classes.  Items 3 and 4 may be less clear.  There are well-known results describing the structure of countable Abelian $p$-groups \cite{kaplansky} in terms of the Ulm sequence and the dimension of the divisible part.  Similarly, there are results of Ketonen \cite{Ket}, \cite{P} describing the structure of countable Boolean algebras.  Theorem \ref{thm3.1} says that the invariants for computable Abelian $p$-groups, and computable Boolean algebras are not not uniformly simple.  For the classes in the following theorem, which is proved in \cite{C1}, our intuition may not say whether there should be a structure theorem, but the proofs make it clear that these are as bad as linear orderings or undirected graphs.  

\begin{thm}
For the following classes $K$, $E(K)$ is
$m$-complete $\Sigma^1_1$.

\begin{enumerate}
\item Ordered real closed fields (not necessarily Archimedean),
\item Arbitrary fields of any fixed characteristic.
\end{enumerate}
\end{thm}

The proof for 1 uses a lemma of van den Dries.  For an ordered real closed field $F$, 
we define an equivalence relation $\sim$, where $a\sim b$ if there are polynomials 
$p(x)$ and $q(x)$, with rational coefficients, such that 
$p(a)\geq b$ and $q(b)\geq a$.  The equivalence classes are intervals, so we have an 
induced ordering on $F/_\sim$.  We can pass from 
a computable linear ordering $L$ to a computable real closed field 
$F(L)$ such that $\{x \in F : x \geq 0\}/\sim$ has order type $L$.  

For 2, we get the characteristic zero case from 1.  An alternate proof of that case,
which extends to arbitrary characteristic, uses the Friedman and Stanley embedding 
of undirected graphs in fields (see Section 4).

\bigskip

Friedman and Stanley conjectured (in 1989) that
the class of countable torsion-free Abelian groups is Borel complete
\cite{F-S}.  The conjecture remains open, although Hjorth \cite{Hj} succeeded
in showing that the isomorphism relation on this class is not Borel.  In \cite{Cthes}, 
Hjorth's proof is adapted to give the following result.  

\begin{thm} 

Let $K$ be the class of torsion-free Abelian groups.  Then
$E(K)$ is not hyperarithmetical. 
\end{thm}

The proof involves coding trees in torsion-free Abelian groups (Hjorth used slightly
different structures).  For each $\alpha$, we can show that $E(K)$ is $\Delta^0_\alpha$ 
hard, but the coding method varies with $\alpha$, so we do not get the 
fact that the isomorphism problem is $m$-complete $\Sigma^1_1$.  

\bigskip

So far, we have considered classes for which the isomorphism problem is complicated.  Now, we turn to classes where the isomorphism problem should be less complicated.  We consider classes with simple invariants.  Perhaps the nicest examples are vector spaces over a fixed field and algebraically closed fields of fixed characteristic.  For these classes, results on the isomorphism problems match our intuition.  For the vector spaces, we assume that the relevant field is infinite, so that the class has infinite members
of finite as well as infinite dimension.  We also suppose that the field is computable.  The following result is in \cite{C1}.          

\begin{thm}
\label{vsacf} 

For the following classes $K$, $E(K)$ is $m$-complete $\Pi^0_3$.
\begin{enumerate}
\item Vector spaces over a fixed infinite computable field,
\item Algebraically closed fields of fixed characteristic.

\end{enumerate}
\end{thm}

\begin{proof}In each case, we show that $E(K)$ is
$\Pi^0_3$ by noting that $I(K)$ is $\Pi^0_2$ and the following relation, on pairs $(a,b)$ in $I(K)$, is $\Pi^0_3$:  
\begin{quotation}\noindent
``For all $n$, the structure
$\mathcal{A}_a$ has at least $n$ independent elements if and only if
$\mathcal{A}_b$ does.''
\end{quotation}  
For completeness, we show that for any $\Pi^0_3$ set $S$, there is a uniformly
computable sequence $(\mathcal{A}_n)_{n\in\omega}$ in $K$ such that
$n\in S$ iff $\mathcal{A}_n$ has infinite dimension. \end{proof}

In \cite{Cthes}, Theorem \ref{vsacf} is generalized to the class $K$ of models of a strongly minimal theory satisfying some extra conditions.  The models are classified by dimension.  We add a model-theoretic condition ($acl(\emptyset)$ is infinite) to guarantee that there are infinite models of finite as well as infinite dimension.  We also need some effectiveness conditions (e.g., $T$ is decidable with effective elimination of quantifiers).         

It is also relatively simple to classify the computable Archimedean real closed fields, since each is determined, up to isomorphism, by the Dedekind cuts that are filled.  The following result is in \cite{C1}.

\begin{thm} 

If $K$ is the class of Archimedean real closed fields,
then $E(K)$ is $m$-complete $\Pi^0_3$.
  
\end{thm}

At this point, we return to the class of Abelian $p$-groups.  Although we have an algebraic classification of countable Abelian $p$-groups, we have seen that the isomorphism problem is $m$-complete $\Sigma^1_1$.  In fact, there is a single computable Abelian 
$p$-group whose inclusion in the class guarantees this.  Below, we consider subclasses of Abelian $p$-groups that are ``reduced'' and have Ulm sequences of bounded length.  For these subclasses, the isomorphism problem is hyperarithmetical, with complexity increasing with the bound on the length.    

Below, we summarize the definitions leading to the Ulm sequence.  For more information, see \cite{kaplansky}.  An Abelian $p$-group is \emph{reduced} if it has no nontrivial divisible subgroup.  Let $G$ be a countable reduced Abelian
$p$-group.  We begin with the sequence of subgroups $G_\alpha$, where $G_0 = G$, $G_{\alpha+1} = pG_\alpha$, and for limit $\alpha$, $G_\alpha = \cap_{\beta<\alpha} G_\beta$.  Let $P$ be the set of elements of $G$ of order $p$, and let $P_\alpha = P\cap G_\alpha$.  Since $G$ is reduced, there is a countable ordinal $\alpha$ such that $G_\alpha = \{0\}$.  The least such $\alpha$ is the \emph{length} $\lambda(G)$.  For each $\beta <\alpha$, 
$P_\beta/P_{\beta+1}$ is a vector space over the field with $p$ elements.  We write $u_\beta(G)$ for the dimension, which is either finite or $\infty$.  The \emph{Ulm sequence} is the sequence $(u_\beta(G))_{\beta < \lambda(G)}$. 

Ulm's Theorem says that a countable reduced Abelian $p$-group is determined, up to isomorphism, by its Ulm sequence.  Previous papers, including \cite{barker} and
\cite{C2} have given formulas, for any $\alpha$ and $k$, of complexity
(depending on $\alpha$), expressing that the $\alpha$th Ulm invariant
is at least $k$.  Let $K_\alpha$ be the class of reduced Abelian $p$-groups of length
at most $\alpha$.  In \cite{C2}, there are results for these classes.  It is shown that $E(K_\omega)$ is
$\Pi^0_3$, $E(K_{\omega\cdot 2})$ is $\Pi^0_5$, and in general, $E(K_{\omega\cdot\alpha})$ is $\Pi^0_{\hat{\alpha}}$, where $\hat{\alpha} =
\mbox{$\sup\limits_{\gamma < \alpha} (2 \gamma + 3)$}$

We have a problem in expressing a meaningful completeness result.  Until now, we have considered classes $K$ for which $I(K)$ was simpler than $E(K)$.  However, $I(K_\omega)$ is $\Pi^0_3$, the same complexity as $E(K_\omega)$ (and, in general, $I(K_{\omega\cdot\alpha})$ is $\Pi^0_{\hat{\alpha}}$, the same as $E(K_{\omega\cdot\alpha})$) .  It is conceivable that when we decide whether $(a,b)\in E(K_\omega)$, the main difficulty is determining whether $a,b\in I(K)$.  The following definition, stated in \cite{C2}, separates the isomorphism problem from the problem of determining
membership in the class.  

\begin{defn} Suppose $A \subseteq B$, and let $\Gamma$ be a
complexity class (e.g.\ $\Pi^0_3$). 
\begin{enumerate}
\item We say that $A$ is \emph{$\Gamma$ within $B$} if 
  there is some $R \in \Gamma$ such that $A = R \cap B$.
\item We say that \emph{$S \leq_m A$ within $B$} if there is a
  computable $f: \omega \to B$ such that for all $n$, $n \in S$ iff $f(n) \in A$.
\item We say that $A$ is \emph{$m$-complete $\Gamma$ within $B$} if $A$ is
  $\Gamma$ within $B$ and for all $S \in \Gamma$, we have $S \leq_m A$
  within $B$.
\end{enumerate} \end{defn}

The relevant example is where $B = I(K)^2$ for some class $K$, and 
$A = E(K)$.  Note that the condition of being $m$-complete $\Gamma$ within $B$ 
is stronger that simply being $m$-complete $\Gamma$, since the ``negative'' 
examples produced by the reduction must be members of $B$.  We will often write
``Within $K$'' to mean ``Within $I(K)^2$.''  All previous theorems
in this section remain true when ``Within $K$'' is added to the
statement.  Using these definitions, we state the result from \cite{C2}, with the
appropriate completeness.  

\begin{thm}
\label{eapg}
Let $K_{\omega\cdot\alpha}$ be the class of reduced Abelian $p$-groups of
length at most $\omega\cdot\alpha$, and let $\hat{\alpha} =
\mbox{$\sup\limits_{\gamma < \alpha} (2 \gamma + 3)$}$.
Then $E(K_{\omega\cdot\alpha})$ is $m$-complete $\Pi^0_{\hat{\alpha}}$ 
within~$K_\alpha$.

\end{thm}

Now we return to torsion-free groups, to consider a simpler
subclass of those.  The torsion-free Abelian groups are, essentially,
subgroups of vector spaces over the rationals.  The \emph{rank} of
such a group is the dimension of the least $\mathbb{Q}$-vector space
in which the group can be embedded.  For groups of finite rank, the
isomorphism problem is simple because there are simple isomorphisms.

We say that a structure $\mathcal{A}$ is \emph{computably categorical} if for any 
$\mathcal{B} \simeq\mathcal{A}$, there is some computable function $f: \mathcal{B} \to
\mathcal{A}$ which is an isomorphism.  (For more about computable categoricity, see 
\cite{A-K}.)  If all computable members of a
class $K$ are computably categorical, then $E(K)$ must be
$\Sigma^0_3$, at least within $K$.  

Any torsion-free Abelian group of finite rank is computably categorical.  Using this, it is 
possible to prove the following \cite{Cthes}.

\begin{thm} 

Let $K$ be either the class of all torsion-free
Abelian groups of finite rank, or the class of torsion-free
Abelian groups of rank~$1$.  Then $E(K)$ is $m$-complete 
$\Sigma^0_3$ within $K$.

\end{thm}

\subsection{Optimal Scott sentences and index sets}

Classification involves describing the structures in a class, up to
isomorphism.  One way to do this is by giving Scott sentences.  If we
consider only computable structures, it is enough to give pseudo-Scott
sentences.  How do we know when we have an optimal Scott sentence, or
pseudo-Scott sentence?  

For inspiration, we go to the setting of descriptive set theory, with arbitrary countable structures and Scott sentences in $L_{\omega_1\omega}$.  As we said in Section 1, a class $K$ closed under automorphism is Borel iff it is the class of models of some $L_{\omega_1\omega}$ sentence.  More is true.  The class $K$ is $\mathbf{\Sigma^0_\alpha}$, or $\mathbf{\Pi^0_\alpha}$, iff it is the class of models of some $\Sigma_\alpha$, or $\Pi_\alpha$ sentence.  If $K(\mathcal{A})$ is the class of copies of a given $\mathcal{A}$, then the form of the optimal Scott sentence matches the complexity of $K$.  D.\ Miller \cite{DM} showed that if $K(\mathcal{A})$ is $\mathbf{\Delta^0_\alpha}$, then it is $d$-$\mathbf{\Sigma^0_\alpha}$.  A.\ Miller \cite{AM} gave examples of structures $\mathcal{A}$ such that $K(\mathcal{A})$ is $\mathbf{\Sigma^0_\alpha}$, $\mathbf{\Pi^0_\alpha}$, or 
$d$-$\mathbf{\Sigma^0_\alpha}$ for various countable ordinals.  He showed that $K(\mathcal{A})$ cannot be properly $\mathbf{\Sigma^0_2}$.          

Now, we return to the computable setting, with computable structures and computable infinitary sentences.  Let $\mathcal{A}$ be a computable structure, and let $K$ be the class of structures for the language of $\mathcal{A}$.  If $\mathcal{A}$ has a computable $\Pi_\alpha$ Scott sentence, or pseudo-Scott sentence, then it is easy to see that $I(\mathcal{A})$ is $\Pi^0_\alpha$, within 
$I(K)$.  Similarly, if $\mathcal{A}$ has a computable $d$-$\Sigma_\alpha$
Scott sentence, or pseudo-Scott (i.e., one of the form $\varphi\ \&\ \neg{\psi}$, where $\varphi$ and $\psi$ are computable $\Sigma_\alpha$), then $I(\mathcal{A})$ is
$d$-$\Sigma^0_\alpha$ within $I(K)$.  

In \cite{C-H-K-M}, there are preliminary results supporting the thesis that for computable structures $\mathcal{A}$ of various familiar kinds, when we have found an optimal
Scott sentence, or pseudo-Scott sentence, then the complexity of $I(\mathcal{A})$ 
matches that of the sentence, and $I(\mathcal{A})$ is $m$-complete in the complexity class.  We begin with vector spaces.

\begin{thm} Let $K$ be the class of vector spaces over
  $\rat$, let $K_f$ consist of the finite dimensional members of $K$, and let
  $\sA$ be a computable member of $K$.  
  
\begin{enumerate}
\item If $\dim \sA = 0$ then $I(\sA)$ is $m$-complete $\Pi^0_1$ within $K_f$.
\item If $\dim \sA = 1$ then $I(\sA)$ is $m$-complete $\Pi^0_2$ within $K_f$.
\item If $1 < \dim \sA < \aleph_0$ then $I(\sA)$ is $m$-complete 
$d$-$\Sigma^0_2$ within $K_f$.
\item If $\dim \sA = \aleph_0$ then $I(\sA)$ is $m$-complete $\Pi^0_3$
  within $K$.
\end{enumerate}
\end{thm}

For 3, the Scott sentence for $\mathcal{A}$, within $K_f$, says that there are at least $n$ linearly independent elements and there are not at least $n+1$.  For 2, we might first write a Scott sentence of the same kind, but this is not optimal.  The optimal sentence says that any two elements are linearly dependent.  

\bigskip

The situation for Archimedean ordered real closed fields
is similar.

\begin{thm} Let $K$ be the class of computable Archimedean real closed
ordered fields, and let $\sA$ be a computable member of $K$.
\begin{enumerate}
\item If the transcendence degree of $\sA$ is $0$, then $I(\sA)$
  is $m$-complete $\Pi^0_2$ within~$K$.
\item If the transcendence degree of $\sA$ is finite but greater than
  $0$, then $I(\sA)$ is $m$-complete $d$-$\Sigma^0_2$ within $K$.
\item If the transcendence degree of $\sA$ is infinite then $I(\sA)$
  is $m$-complete $\Pi^0_3$ within $K$.
\end{enumerate}
\end{thm}

There is a computable $\Pi_2$ sentence $\varphi$ characterizing the Archimedean real closed ordered fields.  For 1, there is a computable $\Pi_2$ Scott sentence, the conjunction of $\varphi$ and a sentence saying that all elements are algebraic.  For 2, fix a transcendence base, of size $n$, say.  There is a computable $\Sigma_2$ Scott sentence---the conjunction of $\varphi$ and a sentence saying that the cuts of the transcendence base are filled, and there do not exist $n+1$ algebraically independent elements.  For 3, there is a computable $\Pi_3$ Scott sentence---the conjunction of $\varphi$ and a sentence saying the appropriate cuts are filled, and every element lies in one of the cuts that should be filled.

\bigskip      

Now, we turn to Abelian $p$-groups.  The proof of Theorem \ref{eapg} does some index set calculations.  In each case, the completeness was witnessed by a single reduced Abelian $p$-group.  We get the following.  

\begin{prop}
 
Let $\mathcal{A}$ be the reduced Abelian $p$-group of length $\omega\cdot\alpha$, with $u_\gamma(\mathcal{A}) = \infty$ for all $\gamma < \alpha$.  Let 
$\hat{\alpha} = sup_{\omega\gamma < \alpha} (2\gamma+3)$.  Then $I(\mathcal{A})$ is $m$-complete $\Pi^0_{\hat{\alpha}}$ within the class of reduced Abelian $p$-groups of length $\leq\omega\cdot\alpha$.   

\end{prop}

A theorem of Khisamiev \cite{Kh1}, \cite{Kh2} gives a description of the Ulm
invariants of computable Abelian $p$-groups of length less than
$\omega^2$.  For each of these groups, we can determine the complexity of
the index set. 

\begin{thm}

Let $K$ be the class of reduced Abelian $p$-groups of length $\omega~M~+~N$
for some $M, N \in \omega$. Let $\mathcal{A} \in K$.

\begin{enumerate}
\item If $\mathcal{A}_{\omega M}$ is minimal for the given length (i.e., it has the
form $\mathbb{Z}_{p^{N}})$, then $I(\mathcal{A})$ is $m$-complete 
$\Pi_{2M+1}^{0}$ within $K$.

\item If $\mathcal{A}_{\omega M}$ is finite but not minimal for the given
length, then $I(\mathcal{A})$ is\linebreak $m$-complete $d$-$\Sigma _{2M+1}^{0}$
within $K$.

\item If there is a unique $k < N$ such that $u_{\omega M + k}(\mathcal{A})
= \infty$, and for all $m < k$, $u_{\omega M + m}(\mathcal{A}) = 0$, then 
$I(\mathcal{A})$ is $m$-complete $\Pi^0_{2M+2}$ within $K$.

\item If there is a unique $k < N$ such that $u_{\omega M + k}(\mathcal{A})
= \infty$ and for some $m < k$ we have $0 < u_{\omega M + m}(\mathcal{A}) <
\infty$, then $I(\mathcal{A})$ is $m$-complete $d$-$\Sigma^0_{2M+2}$ 
within~$K$.

\item If there exist $m<k<N$ such that 
$u_{\omega M+m}(\mathcal{A})=u_{\omega M+k}(\mathcal{A})=\infty$,
then $I(\mathcal{A})$ is $m$-complete $\Pi _{2M+3}^{0}$ within $K$.
\end{enumerate}
\end{thm}

In most cases here, the first Scott sentence that we write down, describing the Ulm sequence in a straightforward way, is not optimal.  To specify that $u_{\omega M + k}(\mathcal{A}) = \infty$, it suffices, via a Ramsey's Theorem argument, to say that for all $r$, the group
$\mathcal{A}_{\omega M}$ has a subgroup of type
$\mathbb{Z}_{p^{k+1}}$.  This gives a simpler Scott sentence.

\bigskip

In \cite{C-H-K-M}, the models of the
original Ehrenfeucht theory $T$ illustrate a further pattern.  Recall that $T$ is the theory of a dense
linear ordering without endpoints, with an infinite increasing sequence of constants.  The three non-isomorphic models are as follows:

\begin{enumerate}

\item  the prime model, in which the sequence has no upper bound, 

\item  a ``middle'' model, in which the sequence has an upper bound but no least upper bound, and 

\item  the saturated model, in which the sequence has a least upper bound.

\end{enumerate}  

\begin{thm}
Let $K$ be the class of models of the original Ehrenfeucht
theory, and let $\mathcal{A}^1$, $\mathcal{A}^2$, and $\mathcal{A}^3$ be
the prime model, the middle model, and the saturated model, respectively.  Then

\begin{enumerate}
    \item $I(\mathcal{A}^1)$ is $m$-complete $\Pi^0_2$ within $K$.
    \item $I(\mathcal{A}^2)$ is $m$-complete $\Sigma^0_3$ within $K$.
    \item $I(\mathcal{A}^3)$ is $m$-complete $\Pi^0_3$ within $K$.
\end{enumerate}
\end{thm}

In related work, Csima, Montalban, and Shore \cite{CMS}
calculated the index sets for Boolean algebras up to elementary
equivalence.  The calculation was based on elementary invariants isolated by Tarski (see Chapter 7 of \cite{Ko}, or \cite{Goncharov}).   

\section{Computable embeddings} 

In \cite{C-C-K-M}, there is a notion of \emph{computable} embedding, inspired by that of 
\emph{Borel} embedding.  The goal was a definition that would allow
meaningful comparisons even of classes of finite structures.  There
are different possible definitions.  The one chosen in \cite{C-C-K-M}
is essentially uniform enumeration reducibility.  Recall that $B$ is
\emph{enumeration reducible} to $A$, or $B\leq_e A$, if there is a
computably enumerable set $\Phi$ of pairs $(\alpha,b)$, where $\alpha$ is a finite set
and $b$ is a number such that $B = \{b:(\exists\alpha\subseteq
A)\,(\alpha,b)\in\Phi\}$.

\bigskip
\noindent
\textbf{Note}:  Given $\Phi$, for each set $A$, there is a unique set $B$ such that $B\leq_e A$ via
$\Phi$.  Thus, $\Phi$ yields a function from $P(\omega)$ to $P(\omega)$. 

\bigskip

As in the setting of descriptive set theory, we consider structures
that need not be computable.  We suppose that each structure has
universe a subset of $\omega$.  For simplicity, we suppose that the
language consists of finitely many relation symbols---we could allow more
general languages.  Each class consists of structures for a single
language, and, modulo the restriction on the universes, each class is
closed under isomorphism.   

\begin{defn}  

Let $K,K'$ be classes of structures.  A \emph{computable
  transformation} from $K$ to $K'$ is a computably enumerable set $\Phi$ of pairs
$(\alpha,\varphi)$, where 
\begin{enumerate}
\item $\alpha$ is a finite set of sentences
appropriate to be in the atomic diagram of a member of $K$, 
\item $\varphi$ is a sentence appropriate to be in the atomic diagram of a
member of $K'$, and
\item for each $\mathcal{A}\in K$, there exists $\mathcal{B}\in K'$
  such that $\Phi(D(\mathcal{A})) = D(\mathcal{B})$; 
i.e., $D(\mathcal{B})\leq_e D(\mathcal{A})$ via $\Phi$.
\end{enumerate}

\end{defn}
  
Identifying the structures with their atomic diagrams, we write $\Phi(\mathcal{A}) = \mathcal{B}$.  

\begin{defn}  

Let $K,K'$ be classes of structures.  A \emph{computable embedding of
  $K$ in $K'$} is a computable transformation $\Phi$ such that for all
$\mathcal{A},\mathcal{A}'\in K$, $\mathcal{A}\cong\mathcal{A}'$
iff $\Phi(\mathcal{A})\cong\Phi(\mathcal{A}')$ 

\end{defn}

We write $K\leq_c K'$ if there is a computable embedding of $K$ into $K'$.

\bigskip
\noindent
\textbf{Example 1}.  Let $FLO$ be the class of finite linear orders, and let $FVS$ be
the class of finite-dimensional vector spaces over $\mathbb{Q}$.  Then $FLO\leq_c FVS$.  

\bigskip
\noindent
\begin{proof}

Let $\mathcal{V}$ be a computable vector space, with a computable sequence $(b_k)_{k\in\omega}$
of linearly independent elements.  For each finite set $S$, let $\mathcal{V}(S)$ be the linear span of $\{b_k:k\in S\}$.  Let $\Phi$ be the set of pairs $(\alpha,\varphi)$ such that for some finite linear order
$\mathcal{L}$, with universe $S$, $\alpha = D(\mathcal{L})$ and
$\varphi\in D(\mathcal{V}(S))$.

\end{proof}

\bigskip
\noindent
\textbf{Example 2}.  Let $UG$ be the class of undirected graphs, and
let $F_\chi$ be the class of fields of any characteristic $\chi$ ($0$, or
$p$, for some prime $p$).  Then $UG\leq_c F_\chi$.

\begin{proof} 

We adapt a Borel embedding given by Friedman and Stanley \cite{F-S}.
Let $\mathcal{F}$ be a computable algebraically closed field of
characteristic $\chi$, with a computable sequence $(b_k)_{k\in\omega}$ of
algebraically independent elements.  For $\chi\not= 2$, for each graph
$\mathcal{G}$, let $\mathcal{F}(\mathcal{G})$ be the subfield of
$\mathcal{F}$ generated by the elements
\[ \left( \bigcup\limits_{k \in \mathcal{G}} acl(b_k)\right) \cup
\{\sqrt{d_i+d_j} |(d_i,d_j) \in R\}\]
where $R$ is the set of all pairs of elements $(d_i, d_j)$ where $i,j$
are elements of $\mathcal{G}$ joined by an edge, and we have $acl(d_i)
= acl(b_i)$ and $acl(d_j) = acl(b_j)$.\footnote{Friedman and Stanley
  put in only $\sqrt{b_i+b_j}$.  This is enough for a computable
  embedding.  The added elements $\sqrt{d_i+d_j}$ give the embedding
  a further property of ``rank preservation,'' which will be helpful
  in section 5.}  For $\chi = 2$, we
use cube roots instead of square roots to indicate edges.  Let $\Phi$
be the set of pairs $(\alpha,\varphi)$, where for some finite
undirected graph $\mathcal{G}$, $\alpha = D(\mathcal{G})$ and
$\varphi\in D(\mathcal{F}(\mathcal{G}))$.  It is not at all easy to
see that if $\Phi(\mathcal{A}) \cong \Phi(\mathcal{A}')$, then
$\mathcal{A} \cong \mathcal{A}'$.  The difficulty is showing that square
roots of the kind representing edges do not appear in the field by accident.  
One of several proofs in the
literature is given in \cite{shapiro}.

\end{proof}

\bigskip
\noindent
\textbf{Example 3}.  Let $LO$ be the class of linear orderings.  Then $UG\leq_c LO$.

\begin{proof}

Again we adapt a Borel embedding given by Friedman and Stanley
\cite{F-S}.  Let $\mathcal{L}$ be $Q^{<\omega}$, under the
lexicographic ordering.  Let $(Q_a)_{a\in\omega}$ be a uniformly
computable disjoint family of dense subsets of $\mathbb{Q}$.  Let
$(t_n)_{n\in\omega}$ be a computable list of the atomic types of
tuples in the language of graphs.  For a graph $\mathcal{G}$, we let
$\mathcal{L}(\mathcal{G})$ be the sub-ordering of $\mathcal{L}$
consisting of the sequences $r_0q_1r_1q_2r_2\ldots q_nr_n k$ such that
for some tuple $(a_1,\ldots,a_n)$ realizing type $t_m$ in
$\mathcal{G}$, $q_i\in Q_{a_i}$, for $i < n$, $r_i\in Q_0$, $r_n\in
Q_n$, and $k$ is a natural number with $k < m+1$.  We let $\Phi$ 
be the set of pairs $(\alpha,\varphi)$, where for some finite undirected graph
$\mathcal{G}$, $\alpha = D(\mathcal{G})$, and $\varphi\in
D(\mathcal{L}(\mathcal{G}))$. 

\end{proof}

The following fact has no analogue in the setting of Borel embeddings.

\begin{prop}
\label{prop4.1}

Suppose $K\leq_c K'$ via $\Phi$.  If $\mathcal{A},\mathcal{B}\in K$, where $\mathcal{A}\subseteq\mathcal{B}$, then $\Phi(\mathcal{A})\subseteq\Phi(\mathcal{B})$.

\end{prop} 

It follows that if $K$ contains a properly increasing chain of structures of length $\alpha$, and $K'$ has no such chain, then $K\not\leq_c K'$.  Using this proposition, we can see that there is no computable embedding of the class $LO$ of linear orderings in the class of vector spaces. 

\begin{prop}

The relation $\leq_c$ is a partial order on the set $\mathcal{C}$ of all classes
of structures.

\end{prop} 

As usual, from the partial ordering $\leq_c$ on $\mathcal{C}$, we get an equivalence relation $\equiv_c$, also on $\mathcal{C}$, where $K\equiv_c K'$ iff 
$K\leq_c K'$ and $K'\leq_c K$.  Moreover, $\leq_c$ induces a partial order on the $\equiv_c$-classes, which we call \emph{$c$-degrees}.  Let $\mathbf{C}$ be the resulting degree structure 
$(\mathcal{C}/_{\equiv_c},\leq_c)$. 

\subsection{Using computable embeddings to compare classes} 

Familiar classes of finite structures seem to lie in one of two $c$-degrees.  Recall that $FLO$ is the class of finite linear orders.  Let $PF$ be the class of finite prime fields.  

\begin{prop}\
\label{prop4.3}

\begin{enumerate}

\item  $PF <_c FLO$,

\item  the $c$-degree of $PF$ also contains the class of finite cyclic graphs,

\item  the $c$-degree of $FLO$, maximum for classes of finite structures, also contains the following classes:  finite undirected graphs, finite cyclic groups, and finite simple groups.  

\end{enumerate}  

\end{prop}

Considering infinite structures, we come to further landmark $c$-degrees.  Recall that $FVS$ is the class of finite dimensional vector spaces over $\mathbb{Q}$.  Let $LO$ be the class of linear orders, infinite as well as finite.  

\begin{prop}\ 
\label{prop4.4}      

\begin{enumerate}

\item  $FLO <_c FVS <_c LO$,  

\item  the $c$-degree of $LO$, maximum over-all, also contains the class of undirected graphs, and the class of fields of characteristic $0$ (or $p$).       

\end{enumerate}

\end{prop}

The class of vector spaces over $\mathbb{Q}$, the class of algebraically closed fields of fixed characteristic, and the class of models of $Th(Z,S)$ (the integers with successor) share some important features.  We do not know whether they represent the same $c$-degree.   

\bigskip

There are natural questions about the degree structure itself.  Is it a linear order?  Is it a lattice?  The first of these questions is answered in \cite{C-C-K-M}, and the second is answered in \cite{K}.                        

\begin{thm}
\label{thm4.5} 

There is a family of $2^{\aleph_0}$ pairwise
incomparable $c$-degrees.  Therefore, $\mathbf{C}$ is not a linear order.

\end{thm}

The proof of Theorem \ref{thm4.5} uses the following notion.  

\begin{defn}  

We say that sets $A,B$ are \emph{bi-immune}\footnote{The term ``bi-immune'' was used earlier in a quite different way in,
for instance, \cite{J}.} if for any partial computable function $f$,

\begin{enumerate}

\item  if $ran(f\upharpoonright A)$ is infinite, then $ran(f\upharpoonright A)\not\subseteq B$, 

\item  if $ran(f\upharpoonright B)$ is infinite, then $ran(f\upharpoonright B)\not\subseteq A$.

\end{enumerate}

\end{defn}

For each natural number $n$, we have a finite prime field of size $p_n$, where $p_n$ is the $n^{th}$ prime.   Below, we write $PF(X)$ for the class of finite prime
fields having size $p_n$ for $n \in X$.

\begin{lem}
\label{lem4.6} 

If $A,B$ are bi-immune, then $PF(A)\not\leq_c PF(B)$ and\\
$PF(B)\not\leq_c PF(A)$.

\end{lem}

Lemma \ref{lem4.6} reduces the proof of Theorem \ref{thm4.5} to the following.        
  
\begin{lem} 
\label{lem4.7}

There is a family of $2^{\aleph_0}$ pairwise bi-immune sets.

\end{lem}

By relativizing the notion, and the lemma, we obtain further families of incomparable 
$c$-degrees in various intervals. 

\bigskip 

The result below says that $\mathbf{C}$ is not a lattice.    

\begin{thm}
\label{thm4.8}  

There is a pair of $c$-degrees with no join and no meet.  

\end{thm}

\begin{proof}

Let $A,B$ be mutually generic.  The forcing conditions are the pairs $(p_1,p_2)$ with
$p_1,p_2\in 2^{<\omega}$ (representing finite partial characteristic functions for $A,B$, respectively).  Take the $c$-degrees of $PF(A)$ and $PF(B)$.  With some effort, it can be shown that this pair of degrees has neither meet nor join.    

\end{proof}

For computable embeddings, the connection between the classes is much tighter 
than for Borel embeddings.  One piece of evidence for this is Proposition \ref{prop4.1}, on preservation of the substructure relation.  As further evidence, we consider transfer of invariants.  Suppose 
$K\leq_c K'$ via $\Phi$.  If there are nice invariants for elements of $K'$, we may use them to describe elements of $K$.  There are two methods for doing this.  The first is trivial.  We describe 
$\mathcal{A}\in K$ as the structure whose $\Phi$-image satisfies the description of $\Phi(\mathcal{A})$.  The second method is more satisfying, as it allows us to describe $\mathcal{A}$ in its own language.  The following result is from \cite{K}.                

\begin{thm}
\label{thm4.9}

Suppose $K\leq_c K'$ via $\Phi$.  There is a partial computable function taking computable infinitary sentences $\varphi$ in the language of $K'$ to computable infinitary sentences $\varphi^*$ in the language of $K$ such that for all $\mathcal{A}$ in $K$, $\Phi(\mathcal{A})\models\varphi$ iff $\mathcal{A}\models\varphi^*$.  Moreover, if $\varphi$ is computable $\Sigma_\alpha$, then so is $\varphi^*$.

\end{thm}

\begin{proof}

For all $\mathcal{A}\in K$, we can build generic copies.  The forcing conditions are
members of the set $\mathcal{F}_\mathcal{A}$ of finite $1-1$ partial functions from
$\omega$ to $\mathcal{A}$.  From a generic copy $\mathcal{A}^*$ of $\mathcal{A}$, we get
$\Phi(\mathcal{A}^*)\in K'$.  The forcing language, which describes $\Phi(\mathcal{A}^*)$, is quite different from the language of $K$.  However, forcing is definable in the language of $K$.  For each sentence $\sigma$ in the forcing language and each tuple $\overline{b}$ in $\omega$, we can find a formula $Force(\overline{x})$ such that $\mathcal{A}\models Force(\overline{a})$ iff we have $p\in\mathcal{F}_\mathcal{A}$ taking $\overline{b}$ to $\overline{a}$ and $p$ forces $\sigma$.  The formulas defining forcing are in the language of $K$.  Moreover, they do not depend on
$\mathcal{A}$.  The sentence $\varphi^*$ says, in all $\mathcal{A}\in K$, 
$(\exists p\in\mathcal{F}_\mathcal{A})\,[\,p\Vdash$``$\Phi(\mathcal{A}^*)\models\varphi$''$\,]$.

\end{proof}   
                                            
\section{Computable structures of high Scott rank}

We return to computable structures.  By Theorem \ref{thm2.6}, these
structures all have Scott rank at most $\omega_1^{CK}+1$.  We consider
examples of computable structures having different Scott ranks.  We
begin with structures having computable ranks.  The exact value is
unimportant---that depends on our choice of definition.

\begin{prop}
\label{prop5.1}  

For all computable structures $\mathcal{A}$ in the following classes, 
$SR(\mathcal{A}) < \omega_1^{CK}$.

\begin{enumerate}

\item  well orderings,

\item  superatomic Boolean algebras,

\item  reduced Abelian $p$-groups.

\end{enumerate} 

\end{prop}
 
For 2, note that a countable superatomic Boolean algebra can be expressed as a finite join of $\alpha$-atoms, for some countable ordinal $\alpha$.  In a computable Boolean algebra $\mathcal{A}$, any element that bounds $\alpha$-atoms for all computable ordinals $\alpha$ can be split into two disjoint elements with this feature.  It follows that if $\mathcal{A}$ bounds $\alpha$-atoms for all computable ordinals $\alpha$, then it has an atomless subalgebra, so it cannot be superatomic.  For 3, we show that for a computable Abelian $p$-group, if the length is not computable, then there is a non-zero divisible element.  

\bigskip  
    
Next, we consider computable structures of Scott rank $\omega_1^{CK}+1$.  We have already discussed one example, the Harrison ordering.  This gives rise to some other structures.  In particular, the \emph{Harrison Boolean algebra} is the interval algebra of the Harrison ordering.  The \emph{Harrison Abelian $p$-group} is the Abelian $p$-group of length $\omega_1^{CK}$, with all infinite Ulm invariants, and with a divisible part of infinite dimension.  Clearly, the Harrison Boolean algebra has a computable copy.  The Harrison Abelian $p$-group does as well.        

\begin{thm}

The Harrison Boolean algebra, and the Harrison Abelian\linebreak $p$-groups have Scott rank $\omega_1^{CK}+1$.  

\end{thm}

In the Harrison Boolean algebra, any non-superatomic element has Scott rank $\omega_1^{CK}$.  
In the Harrison Abelian $p$-group, any divisible element has Scott rank $\omega_1^{CK}$.  See 
\cite{G-H-K-S} for more about these structures.    

\bigskip

Finally, we turn to Scott rank $\omega_1^{CK}$.  The following is in \cite{M}.

\begin{thm} [Makkai] 
\label{thm5.3} 

There is an arithmetical structure of Scott rank $\omega_1^{CK}$. 

\end{thm}

In \cite{K-Y}, Makkai's example is made computable.

\begin{thm} 
\label{thm5.4}

There is a computable structure of rank $\omega_1^{CK}$. 

\end{thm}

There are two proofs in \cite{K-Y}.  The first simply codes Makkai's example in a computable structure, without examining it.  The second proof re-works Makkai's construction, incorporating suggestions of Shelah and Sacks.  The structures are ``group trees''.  To obtain a group tree, we start with a tree.  Next, we define a family of groups, one for each level of the tree.  Using the groups, we get a new, more homogeneous tree structure.  We then throw away most of the structure.  We retain a family of unary functions, one for each element.  Morozov \cite{Mo} gives a nice description of the construction.  He was interested in the fact that if we start with a computable tree $T$ having a path but no hyperarithmetical path, then the resulting group tree (Morozov calls it a ``polygon'') has non-trivial automorphisms but no hyperarithmetical ones.  The Harrison ordering (as originally constructed) has no hyperarithmetical infinite decreasing sequence, so it also has no non-trivial hyperarithmetical automorphisms.          

\subsection{Trees}

In \cite{C-K-Y}, there are simpler examples of computable structures of rank $\omega_1^{CK}$, which are just trees.  We consider subtrees of $\omega^{<\omega}$ and their isomorphic copies.  We write $\emptyset$ for the top node.  To describe the examples, we need some terminology.  We define \emph{tree rank} for $a\in T$, and for $T$ itself. 

\begin{defn}\   

\begin{enumerate}

\item  $rk(a) = 0$ if $a$ is terminal, 

\item  for an ordinal $\alpha > 0$, $rk(a) = \alpha$ if all successors
  of $a$ have ordinal rank, and $\alpha$ is the first ordinal greater
  than any of these ranks, 

\item  $rk(a) = \infty$ if $a$ does not have ordinal rank.

\end{enumerate}
We let $rk(T) = rk(\emptyset)$.

\end{defn}

\noindent
\textbf{Fact}.  $rk(a) = \infty$ iff there is a path through $a$.   

\bigskip

For a tree $T$, let $T_n$ be the set of elements at level $n$.

\begin{defn}  

The tree $T$ is \emph{rank-homogeneous} if 

\begin{enumerate}

\item  for all $a\in T_n$ and all ordinals $\alpha$, if there exists $b\in T_{n+1}$ with\\ $rk(b) = \alpha < rk(a)$, then $a$ has infinitely many successors $a'$ with\\ $rk(a') = \alpha$. 

\item  for all $a\in T_n$ with $rk(a) = \infty$, $a$ has infinitely many successors $a'$ with $rk(a') = \infty$. 

\end{enumerate}
\end{defn}

The following is clear.  

\begin{prop}
\label{prop5.5}

If $T$ and $T'$ are rank-homogeneous trees, and for all $n$, $T_n$ and $T'_n$ 
represent the same ranks, then $T\cong T'$. 

\end{prop}  

\begin{defn}  

A tree $T$ is \emph{thin} if for all $n$, the set of ordinal ranks of elements of $T_n$ has order type at most $\omega\cdot n$.

\end{defn}

We use thinness in the following way.

\begin{prop}
\label{prop5.6}

If $T$ is a computable thin tree, then for each $n$, there is some computable $\alpha_n$ such that for all $a\in T_n$, if $rk(a) \geq\alpha_n$, then $rk(a) = \infty$. 

\end{prop}   

To show that there is a computable tree of Scott rank $\omega_1^{CK}$, we show the following.

\begin{thm}\
\label{thm5.7}

\begin{enumerate}

\item  There is a computable, thin, rank-homogeneous tree $T$ such that $T$ has a path but no hyperarithmetical path.

\item  If $T$ is a computable, thin, rank-homogeneous tree such that 
$T$ has a path but no hyperarithmetical path, then $SR(T) = \omega_1^{CK}$.

\end{enumerate}

\end{thm}

\begin{proof}

Part 1 is proved using Barwise-Kreisel Compactness.  There is a
$\Pi^1_1$ set of computable infinitary sentences describing the
desired tree, and we show, in several steps, that every $\Delta^1_1$
subset has a computable model.

For Part 2, first we note that the orbit of each tuple in $T$ is defined by a computable infinitary 
formula---this is where we use thinness.  It follows that
$SR(T)\leq\omega_1^{CK}$.  To show that $SR(T)\geq\omega_1^{CK}$, we
show that if $a$ and $b$ are elements of $T$ at the same level, and
are such that $rk(a) = \infty$ and $b$ has
computable tree rank at least $\geq\omega\cdot\alpha$, then $a\equiv^\alpha b$.
In fact, two tuples are $\alpha$-equivalent if they generate
isomorphic subtrees and for each pair of corresponding elements in the
two subtrees, either the tree ranks match, or else both are at least $\omega\cdot\alpha$.

\end{proof}

In \cite{C-G-K}, examples are given of structures of Scott rank $\omega_1^{CK}$ in further familiar classes.  

\begin{thm} 
\label{thm5.8} 

Each of the following classes contains computable structures of Scott rank $\omega_1^{CK}$:

\begin{enumerate}

\item  undirected graphs,

\item  fields of any fixed characteristic,

\item  linear orderings.

\end{enumerate} 

\end{thm}

To prove Theorem \ref{thm5.8}, we use some special computable
embeddings.  

\begin{defn}

Let $\Phi$ be an embedding of $K$ in $K'$.  We say that
$\Phi$ has the \emph{rank preservation property} if for all computable
$\mathcal{A}\in K$, either $\mathcal{A}$ and $\Phi(\mathcal{A})$ both
have computable Scott rank, or else they have the same non-computable
Scott rank.

\end{defn}

\begin{thm}
\label{thm5.9}

For each of the following pairs of classes $(K,K')$, there is a computable embedding
of $K$ in $K'$ with the rank preservation property:  

\begin{enumerate}

\item  $K$ is the class of trees and $K'$ is the class of undirected graphs,

\item  $K$ is the class of undirected graphs and $K'$ is the class of fields, of any desired characteristic,

\item  $K$ is the class of undirected graphs and $K'$ is the class of linear orderings.

\end{enumerate}

\end{thm}

\begin{proof}

For 1, if $\mathcal{G}$ is the graph corresponding to a tree $T$, then we 
have a copy $T_\mathcal{G}$ of $T$, which is definable in $\mathcal{G}$ using 
existential formulas.  Each automorphism of $T_\mathcal{G}$ extends to an 
automorphism of $\mathcal{G}$.  This implies
that if $\mathcal{G}$ has Scott rank less than $\omega_1^{CK}$, or
at most $\omega_1^{CK}$, then so does $T$.  In addition, for each
tuple $\overline{b}$ in $\mathcal{G}$, there is an existential formula defining the orbit
of $\overline{b}$ under automorphisms that fix the elements of $T_\mathcal{G}$.  
This implies that if $T$ has Scott rank less than $\omega_1^{CK}$, or at 
most $\omega_1^{CK}$, then so does $\mathcal{G}$.

For 2, we use the embedding of graphs in fields described in Section 4
(a variant of the one due to Friedman and Stanley).  There is a copy
of the graph sitting in the corresponding field as a quotient of a
definable structure by a definable congruence relation, and we obtain 
rank preservation 
by proving properties analogous to the ones above.

For 3, we use the embedding of graphs in linear orderings described in
Section 4 (due to Friedman and Stanley).  There is no copy of the
graph definable, even as a quotient structure, in the corresponding
linear ordering.  Even so, the orbit of each tuple in the graph
corresponds to an orbit of an element in the linear ordering, and each
orbit in the linear ordering is determined by an existential formula
together with a finite tuple of orbits from the graph.

\end{proof}                
 
\begin{defn}  

A computable structure $\mathcal{A}$ of non-computable rank is \emph{computably approximable} if every computable infinitary sentence $\sigma$ true in $\mathcal{A}$ is also true in some computable $\mathcal{B}\not\cong\mathcal{A}$. 

\end{defn}

For all we know, every computable structure of high rank is computably approximable.  Failure to be computably approximable is equivalent to existence of a computable infinitary pseudo-Scott sentence.           

\bigskip

Sometimes we want a stronger computable approximability.

\begin{defn}  

Suppose $\mathcal{A}$ is computable, with non-computable Scott rank.  We say that $\mathcal{A}$ is \emph{strongly computably approximable} if for all $\Sigma^1_1$ sets $S$, there is a uniformly computable sequence $(\mathcal{C}_n)_{n\in\omega}$ such that if $n\in S$, then $\mathcal{C}_n\cong\mathcal{A}$, and if $n\notin S$, then $\mathcal{C}_n$ has computable rank.  

\end{defn}    

\bigskip
\noindent
\textbf{Example}:  The Harrison ordering is strongly computably approximable.

\begin{thm}
\label{thm5.10}

There is a computable tree $T$ such that $SR(T) = \omega_1^{CK}$ and $T$ is strongly computably approximable.

\end{thm}  

The proof of Theorem \ref{thm5.10} in \cite{C-K-Y} uses the technical proposition below.  The proposition below refers to Kleene's system for assigning numbers to ordinals.   The system $\mathcal{O}$, with a partial ordering $<_\mathcal{O}$, assigns many notations to a single ordinal.  However, for path $P$ through $\mathcal{O}$, there is a unique notation $a$ for each computable ordinal $\alpha$.  In the proof of Theorem \ref{thm5.10}, the path $P$ comes from an initial segment of a Harrison ordering.  We write $|a|$ for the ordinal with notation $a$.     

\begin{prop}
\label{thm5.11}  

Let $P$ be a path through $\mathcal{O}$.  Let $\mathcal{A}$ be a
computable structure.  Let $(\mathcal{A}_a)_{a\in P}$ be a uniformly
computable family of structures.  Suppose that the following are
satisfied:
\begin{enumerate}
\item $\mathcal{A}_a$ satisfies the computable $\Sigma_{|a|}$ sentences true in
$\mathcal{A}$.  
\item $(\leq_a)_{a\in P}$ are reflexive, transitive
binary relations on tuples from the structures $\mathcal{A}_a$
\item the relations have the \emph{back-and-forth property}; i.e.,
$\leq_0$ preserves satisfaction of quantifier-free formulas, and if
$(\mathcal{A}_a,\overline{a})\leq_b (\mathcal{A}_a',\overline{a}')$, and 
 $b' <_\mathcal{O} b$, then for all $\overline{c}'$, there exists $\overline{c}$ such that 
 $(\mathcal{A}_{a'},\overline{a}',\overline{c}')\leq_{b'}
(\mathcal{A}_a,\overline{a},\overline{c})$, and 
\item the restriction of $\leq_b$ to pairs $(\overline{a},\overline{a}')$, where
$\overline{a}\in\mathcal{A}_a$ and $\overline{a}'\in\mathcal{A}_{a'}$, is c.e.\ uniformly in $b,a,a'$.
\end{enumerate}
Then for any $\Sigma^1_1$ set $S$, there is a uniformly computable
sequence $(\mathcal{C}_n)_{n\in\omega}$ such that if $n\in S$, then
$\mathcal{C}_n\cong\mathcal{A}$, and if $n\notin S$, then
$\mathcal{C}_n\cong\mathcal{A}_a$, for some $a\in P$.
 
\end{prop} 
 
\begin{proof}

The proof of Proposition \ref{thm5.11} is a priority construction with actions described by a $\Pi^1_1$ set $\Gamma$ of computable infinitary sentences.  To show that every $\Delta^1_1$ subset of $\Gamma$ is satisfied, we use an infinitely nested priority construction.  We then apply Barwise-Kreisel Compactness.

\end{proof}      

The conditions of Proposition \ref{thm5.11} seem very strong.  However, with effort, we can show that they are satisfied by a special tree of Scott rank $\omega_1^{CK}$, and a family of approximating trees of computable rank.  

\bigskip
  
In \cite{C-G-K}, Theorem \ref{thm5.10} is combined with Theorem \ref{thm5.9}, to produce further structures of Scott rank $\omega_1^{CK}$ that are strongly computably approximable.     

\begin{thm}
\label{thm5.12}  

There are strongly computably approximable structures of Scott rank $\omega_1^{CK}$
 in each of the following classes:  

\begin{enumerate}

\item  undirected graphs, 

\item  fields of any characteristic,

\item  linear orderings. 

\end{enumerate}
Moreover, in each case, the approximating structures have computable Scott rank.

\end{thm}

\end{document}